\documentclass{article}
\usepackage{amssymb,latexsym}
\usepackage{amsmath}
\usepackage{graphics}
\usepackage{graphicx}
\newtheorem{theorem}{Theorem}

\newtheorem{definition}{Definition}

\newenvironment{proof}[1][Proof]{\noindent\textbf{#1.} }{\ \rule{0.5em}{0.5em}}
\numberwithin{equation}{section} \numberwithin{theorem}{section}
\numberwithin{corollary}{section}
\numberwithin{definition}{section} \numberwithin{result}{section}

\begin{document}
\title{Involute-Evolute Curves in Galilean Space $G_3$}
\author{A. Z. AZAK, M. AKY\.{I}\~{G}\.{I}T, S. ERSOY \\
 apirdal@sakarya.edu.tr, maky\.{i}g\.{i}t@sakarya.edu.tr, sersoy@sakarya.edu.tr \\
Department of Mathematics, Faculty of Arts and Sciences,\\
Sakarya University, Sakarya, 54187 TURKEY}

\maketitle
\begin{abstract}
In this paper, definition of involute-evolute curve couple in
Galilean space is given and some well-known theorems for the
involute-evolute curves are obtained in 3-dimensional Galilean
space.


\end{abstract}

\section{Introduction}\label{S:intro}
C. Huggens discovered involutes while trying to build a more
accurate clock, \cite{Boy}. Later, the relations Frenet apparatus
of involute-evolute curve couple in the space $\mathbb{E}^3$ were
given in \cite{H.H.H.}. A. Turgut examined involute-evolute curve
couple in $\mathbb{R}^n$, \cite{Tur}. At the beginning of the
twentieth century, Cayley-Klein discussed Galilean geometry which
is one of the nine geometries of projective space. After that, the
studies with related to the curvature theory were maintained
\cite{Ros, Pav, Sps} and the properties of the curves in the
Galilean space were studied in \cite{ogr}. In this paper, we
define involute-evolute curves couple and give some theorems and
conclusions, which are known from the classical differential
geometry, in the three dimensional Galilean space $G_3$. We hope
these results will be helpful to mathematicians who are
specialized on mathematical modeling.

\section{Preliminaries}\label{Sec2}
The Galilean space $G_3$ is a Cayley-Klein space equipped with the
projective metric of signature (0,0,+,+), as in \cite{Mol}. The
absolute figure of the Galilean Geometry consist of an ordered
triple $\{ w,f,I\}$, where $w$ is the ideal (absolute) plane, $f$
is the line (absolute line) in $w$ and $I$ is the fixed elliptic
involution of points of $f$, \cite{Sps}. In the non-homogeneous
coordinates the similarity group $H_8$ has the form
\begin{equation}\label{2.1}
\begin{array}{l}
 \overline x  = a_{11}  + a_{12} x \\
 \overline y  = a_{21}  + a_{22} x + a_{23} y\cos \varphi  + a_{23} z\sin \varphi  \\
 \overline z  = a_{31}  + a_{32} x - a_{23} y\sin \varphi  + a_{23} z\cos \varphi  \\
 \end{array}
\end{equation}
where $a_{ij}$ and $\varphi$ are real numbers, \cite{Pav}.\\
In what follows the coefficients $a_{12}$ and $a_{23}$ will play
the special role. In particular, for $a_{12}  = a_{23}  = 1$,
(\ref{2.1}) defines the group $B_6  \subset H_8$ of isometries of
Galilean space $G_3$.\\
In $G_3$ there are four classes of lines:\\
\textit{i)} (proper) non-isotropic lines- they don't meet the
absolute line $f$.\\
 \textit{ii)} (proper) isotropic lines- lines that don't belong to the plane $w$ but meet the absolute line $f$.\\
 \textit{iii)} unproper non-isotropic lines-all lines of $w$ but $f$.\\
 \textit{iv)} the absolute line $f$.\\
Planes $x =$constant are Euclidean and so is the plane $w$. Other
planes are isotropic.\\
If a curve $C$ of the class $C^r \,(r \ge 3)$ is given by the
parametrization
\[
r = r(x,y(x),z(x))
\]
then $x$ is a Galilean invariant the arc length on $C$.\\
The curvature is
\[
\kappa  = \sqrt {y''(x)^2  + z''(x)^2 }
\]
and torsion
\[
\tau  = \frac{1}{{\kappa ^2 }}\det (r'(x),r''(x),r'''(x)).
\]
The orthonormal trihedron is defined
\[
\begin{array}{l}
 T(s) = \alpha '(s) = (1,y'(s),z'(s)) \\
 N(s) = \frac{1}{{\kappa  (s)}}(0,y''(s),z''(s)) \\
 B(s) = \frac{1}{{\kappa  (s)}}(0, - z''(s),y''(s)). \\
 \end{array}
\]
The vectors $T, N, B$ are called the vectors of tangent, principal
normal and binormal line of $\alpha$, respectively. For their
derivatives the following Frenet formulas hold
\[
\begin{array}{l}
 T'(s) = \kappa(s)N(s) \\
 N'(s) = \tau(s)B(s) \\
 B'(s) =  - \tau(s)N(s), \\
 \end{array}
\]
\cite{Kam}.\\
Galilean scalar product can be written as
\[
\left\langle {u_1 ,u_2 } \right\rangle  = \left\{ \begin{array}{l}
 \,\,\,\,\,\,\,x_1 x_2 \;\;\,\;\;\;\;,\;\;\;if\;\;x_1  \ne 0\; \vee \;x_2  \ne 0 \\
 y_1 y_2  + z_1 z_2 \;,\;\;if\;\;x_1  = 0\; \wedge \;x_2  = 0 \\
 \end{array} \right.
\]
where $u_1  = (x_1 ,y_1 ,z_1 )$ and $u_2  = (x_2 ,y_2 ,z_2 )$. It
leaves invariant the Galilean norm of the vector $u = (x,y,z)$
defined by
\[\left\| u \right\| = \left\{ \begin{array}{l}
 \,\,\,\,\,\,\,\,\,\,\,x\;\;\;\;\;\;\;\;\;,\;x \ne 0 \\
 \sqrt {{y^2} + {z^2}} \;\;\,,\;x = 0. \\
 \end{array} \right.\]
Let $\alpha$ be a curve given first by

\begin{equation}\label{2.2}
\begin{array}{l}
 \alpha :I \to G_3 ,\;\;\;\;\;\;I \subset \mathbb{R} \\
 \;\;\;\;\;t \to \alpha (t) = (x(t),y(t),z(t)) \\
 \end{array}
\end{equation}
where $x(t),y(t),z(t) \in C^3$ (the set of three-times
continuously differentiable functions) and $t$ run through a real
interval.

\section{Involute- Evolute Curves in Galilean Space}\label{Sec3}
In this section, we give a definition of involute-evolute curve
and obtain some theorems about these curves in $G_3$.
\begin{definition}\label{T3.1.}
Let $\alpha$ and $\alpha ^*$ be two curves in the Galilean space
$G_3$. The curve $\alpha ^*$ is called involute of the curve
$\alpha$ if the tangent vector of the curve $\alpha$ at the point
$\alpha (s)$ passes through the tangent vector of the curve
$\alpha ^*$ at the point $\alpha ^* (s)$ and
\[
\left\langle {T,T^* } \right\rangle  = 0,
\]
where $\left\{ {T,N,B} \right\}$ and $\left\{ {T^* ,N^* ,B^* }
\right\}$ are Frenet frames of $\alpha$ and $\alpha ^*$,
respectively. Also, the curve $\alpha$ is called the evolute of
the curve $\alpha ^*$.
\end{definition}
This definition suffices to the define this curve mate as
\[
\alpha ^*  = \alpha  + \lambda T
\]
(see Figure 1).
\begin{center}
\hfil\scalebox{1}{\includegraphics{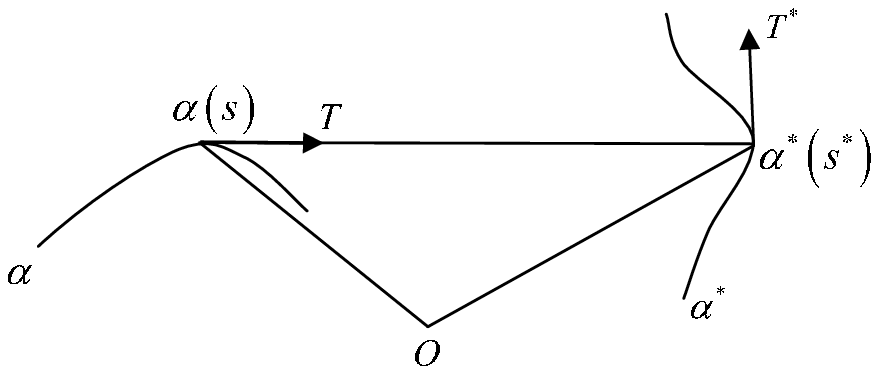}}\hfil\\
\scriptsize{Figure 1. Involute-evolute curves}\\
\end{center}
\begin{theorem}\label{T3.1.}
Let $\alpha$ and $\alpha ^*$ be two curves in the Galilean space
$G_3$. If the curve $\alpha ^*$ is an involute of the curve
$\alpha$, then the distance between the curves $\alpha$ and
$\alpha ^*$ is constant.
\end {theorem}

\begin{proof}
From definition of involute-evolute curve couple, we know
\begin{equation}\label{3.1}
\alpha ^* (s) = \alpha (s) + \lambda (s)T(s)
\end{equation}
Differentiating both sides of the equation (\ref{3.1}) with
respect to $s$ and use the Frenet formulas, we obtain
\[
T^* (s) = T(s) + \frac{{d\lambda }}{{ds}}T(s) + \lambda (s)\kappa
(s)N(s).
\]
Since the curve $\alpha ^*$ is involute of $\alpha$, $\left\langle
{T,T^* } \right\rangle  = 0.$\\
Then we have
\begin{equation}\label{3.2}
\frac{{d\lambda }}{{ds}} + 1 = 0.
\end{equation}
From the last equation, we easily get
\begin{equation}\label{3.3}
\lambda (s) = c - s
\end{equation}
where $c$ is constant. Thus, the equation (\ref{3.1}) can be
written as
\begin{equation}\label{3.4}
\alpha ^* (s) - \alpha (s) = (c - s)T(s).
\end{equation}
Taking the norm of the equation (\ref{3.4}), we reach
\begin{equation}\label{3.5}
\left\| {\alpha ^* (s) - \alpha (s)} \right\| = \left| {c - s}
\right|.
\end{equation}
This completes the proof.
\end{proof}\\

\begin{theorem}\label{T3.2.}
Let $\alpha$ and $\alpha ^*$ be two curves in Galilean space
$G_3$. $\kappa ,\tau$ and $\kappa ^* ,\tau ^*$ be the curvature
functions of $\alpha$ and $\alpha ^*$, respectively. If $\alpha$
is evolute of $\alpha ^* (s)$ then there is a relationship
\[
\kappa ^*  = \frac{\tau }{{(c - s)\kappa }},
\]
where $c$ is constant and $s$ is arc length parameter of $\alpha$.
\end {theorem}

\begin{proof}
Let Frenet frames be $\left\{ {T,N,B} \right\}$ and $\left\{ {T^*
,N^* ,B^* } \right\}$ at the points $\alpha (s)$ and $\alpha ^*
(s)$, respectively. Differentiating both sides of equation
(\ref{3.4}) with respect to $s$ and using Frenet formulas, we have
following equation
\begin{equation}\label{3.6}
T^* (s)\frac{{ds^* }}{{ds}} = (c - s)\kappa (s)N(s)
\end{equation}
where $s$ and $s^*$ are the arc length parameter of the curves
$\alpha$ and $\alpha ^*$, respectively. Taking the norm of the
equation (\ref{3.6}), we reach
\begin{equation}\label{3.7}
T^*  = N
\end{equation}
and
\begin{equation}\label{3.8}
\frac{{ds^* }}{{ds}} = (c - s)\kappa (s).
\end{equation}
By taking the derivative of equation (\ref{3.7}) and using the
Frenet formulas and equation (\ref{3.8}), we obtain
\begin{equation}\label{3.9}
\kappa ^* N^*  = \frac{\tau }{{(c - s)\kappa }}B.
\end{equation}
From the last equation, we get
\[
\kappa ^*  = \frac{\tau }{{(c - s)\kappa }}.
\]
\end{proof}\\

\begin{theorem}\label{T3.3.}
Let $\alpha$ be the non-planar evolute of curve $\alpha ^*$, then
$\alpha$ is a helix.
\end {theorem}

\begin{proof}
Under assumption $s$ and $s^*$ are arc length parameter of the
curves $\alpha$ and $\alpha ^*$, respectively. We take the
derivative of the following equation with respect to $s$
\[
\alpha ^* (s) = \alpha (s) + \lambda (s)T(s)
\]
we obtain that
\[
T^* \frac{{ds^* }}{{ds}} = \lambda \kappa N.
\]
$\left\{ {T^* ,N} \right\}$ are linearly dependent. We may define
function as
\[
f = \left\langle {T,T^*  \wedge N^* } \right\rangle
\]
and take the derivative of the function $f$ with respect to $s$,
we obtain
\begin{equation}\label{3.10}
f' =  - \tau \left\langle {T,{N^*}} \right\rangle.
\end{equation}
From the equation (\ref{3.10}) and the scalar product in Galilean
space, we have
\[
f' = 0.
\]
That is,
\[
f = const.
\]
The velocity vector of the curve $\alpha$ always composes a
constant angle with the normal of the plane which consist of
$\alpha ^*$. Then the non-planar evolute of the curve $\alpha ^*$
is a helix.
\end{proof}\\

\begin{theorem}\label{T3.4.}
Let the curves $\beta$ and $\gamma$ be two evolutes of $\alpha$ in
the Galilean space $G_3$. If the points $P_1$ and $P_2$ correspond
to the point of $\alpha$, then the angle $\mathop {P_1 PP_2
}\limits^ \wedge$ is constant.
\end {theorem}

\begin{proof}
Let's assume that the curves $\beta$ and $\gamma$ be two evolutes
of $\alpha$(see Figure 2). And let the Frenet vectors of the
curves $\alpha ,\,\,\beta$ and $\gamma$ be $\left\{ {T,N,B}
\right\}$, $\left\{ {T^* ,N^* ,B^* } \right\}$ and $\left\{
{\tilde T,\tilde N,\tilde B} \right\}$, respectively.\\
\begin{center}
\hfil\scalebox{1}{\includegraphics{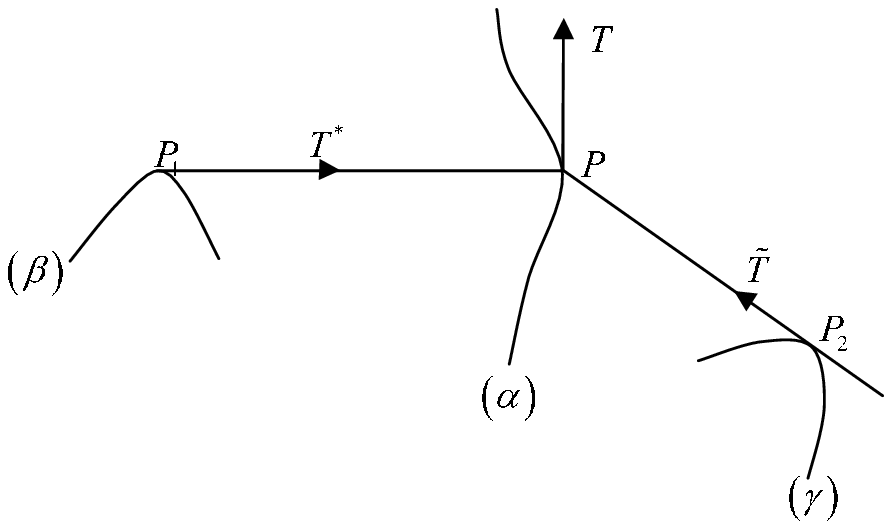}}\hfil\\
\scriptsize{Figure 2. Evolutes of $\alpha$ curve}\\
\end{center}
Following the same way in the proof of the Theorem 3.3, it is
easily seen that $\left\{ {T,N^* } \right\}$ and $\left\{
{T,\tilde N} \right\}$ are linearly
dependent.\\
Thus,
\begin{equation}\label{3.11}
\left\langle {T,T^* } \right\rangle  = 0
\end{equation}
and
\begin{equation}\label{3.12}
\left\langle {T,\tilde T} \right\rangle = 0.\
\end{equation}
We define a function $f$ as an angle between tangent vector $T^*$
and $\tilde T$, that is,
\begin{equation}\label{3.13}
f(s) = \left\langle {{T^*},\tilde T} \right\rangle
\end{equation}
Then, differentiating equation (\ref{3.13}) with respect to $s$,
we have
\begin{equation}\label{3.14}
f'(s) = {\kappa ^*}\frac{{d{s^*}}}{{ds}}\left\langle {{N^*},\tilde
T} \right\rangle + \tilde \kappa \frac{{d\tilde
s}}{{ds}}\left\langle {{T^*},\tilde N} \right\rangle
\end{equation}
where $s$, $s^*$ and $\tilde s$ are arc length parameter of the
curves $\alpha ,\,\,\beta$ and $\gamma$, respectively. Also
$\kappa ,\,\,\kappa ^*$ and $\widetilde\kappa$ are the curvatures
of the curves $\alpha ,\,\,\beta$ and $\gamma$, respectively.
Considering the equations (\ref{3.11}), (\ref{3.12}) and
(\ref{3.14}), we get
\[
f'(s) = 0.
\]
This means that
\[
f = const.
\]
So, the proof is completed.
\end{proof}\\

\end {document}